\documentclass[12pt]{article}
\usepackage{mathrsfs}
\usepackage{mathrsfs}
\usepackage{amsmath}

\def\be{\begin{equation}} \def\ee{\end{equation}}

\def\MCH#1#2{\setbox0=\hbox{\raise#1\hbox{#2}}\smash{\box0}}

\begin{document}

\baselineskip=16pt
\title{\bf The Bahadur representation for sample quantiles under dependent sequence}
\author{Wenzhi Yang ~~~Shuhe Hu\thanks{Corresponding author. E-mail: hushuhe@263.net (S.H. Hu); wzyang827@163.com (W.Z. Yang);
wxjahdx2000@126.com (X.J. Wang)}~~~~Xuejun Wang\\
 {\small School of Mathematical Science, Anhui University, Hefei, Anhui 230039, China}}
\maketitle

\begin{minipage}{155mm}
\indent{\bf Abstract:}~~On the one hand, we investigate the Bahadur representation for sample
quantiles under $\varphi$-mixing sequence with $\varphi
(n)=O(n^{-3})$ and obtain a rate as
$O(n^{-\frac{3}{4}}\log n)$, $a.s.$. On the other hand, by relaxing
the condition of mixing coefficients to $\sum\nolimits_{n=1}^\infty\varphi^{1/2}
(n)<\infty$, a rate $O(n^{-1/2}(\log n)^{1/2})$, $a.s.$, is also
obtained.\\
{\bf Keywords:}~Bahadur representation; sample quantiles; mixing sequence\\
{\it AMS 2000 Mathematical Subject Classification:}~~62F12
\end{minipage}

 \setcounter{equation}{0}
\section{Introduction}

\hspace*{0.5cm} Assume that $\{X_n,n\geq1\}$ is a sequence of random
variables defined on a fixed probability space $(\Omega,
\mathcal{F}, P)$ with a common marginal distribution function
$F(x)=P(X_1\leq x)$. Let $F$ be a distribution function (continuous
from the right, as usual). For $0<p<1$, the $p$th quantile of $F$ is
defined as
$$\xi_p=\inf\{x:F(x)\geq p\}$$
and is alternately denoted by $F^{-1}(p)$. The function $F^{-1}(t)$,
$0<t<1$, is called the inverse function of $F$.

For a sample $X_1,X_2,\cdots,X_n$, $n\geq 1$, let $F_n$ represent
the empirical distribution function based on $X_1,X_2,\cdots,X_n$,
which is defined as $F_n(x)=\frac{1}{n}\sum\nolimits_{i=1}^n
I(X_i\leq x)$, $x\in R$. Here $I(A)$ denotes the indicator function
of the set $A$ and $R$ is the real line. For $0<p<1$, we define
$F_n^{-1}(p)=\inf\{x:F_n(x)\geq p\}$ as the sample $p$th quantile.

Let $\xi_{p,n}=F_n^{-1}(p)$. Bahadur \cite{bahadur 1966} firstly
introduced an elegant representation for sample quantile in terms of
empirical distribution function based on independent and identically
distributed ($i.i.d.$) random variables and obtained the following
result (or see Serfling \cite[Theorem 2.5.1]{serfling 1980})

\textbf{Theorem 1.1.} {\it Let $0<p<1$ and and $\{X_n,n\geq 1\}$ be
a sequence of $i.i.d.$ random variables. Suppose that $F$ is twice
differentiable at $\xi_p$, with $F^{\prime}(\xi_p)=f(\xi_p)>0$. Then
with probability 1,
$$\xi_{p,n}=\xi_p-\frac{F_n(\xi_p)-p}{f(\xi_p)}+O\left(n^{-\frac{3}{4}}(\log
n)^{3/4}\right),~n\rightarrow\infty.$$}

At present, many researchers have extended Bahadur representation
for $i.i.d.$ random variables to the dependent cases. One can refer
to Sen \cite{sen 1972} and Babu and Singh \cite{babu 1978} for
$\varphi$-mixing sequence, Yoshihara \cite{yoshihara 1995} for
$\alpha$-mixing sequence, Zhou and Zhu \cite{zhou 1996} for the
smooth quantile estimator, Sun \cite{sun 2006} for $\alpha$-mixing
sequence, Cheng and Gooijer \cite{cheng 2009} for $M$-estimator
under $\alpha$-mixing sequence, etc. Ling \cite{ling 2008} extended
the results of Sun \cite{sun 2006} to the case of NA sequence and
obtained a rate $O(\tau_n)$ $a.s.$, where $\tau_n\rightarrow 0$ and
$\sqrt{n}\tau_n^2/\log n\rightarrow 0$ as $n\rightarrow\infty$. Li
et al. \cite{li 2011} extended the results of Ling \cite{ling 2008}
to the case of NOD sequence, which is weaker than NA sequence, and
they got a better rate $O(n^{-1/2}(\log n)^{1/2})$ $a.s$. For more
works on Bahadur representation, one can refer to \cite{wang 2011,wang 2016,wu 2005,yang 2014,zhang 2013}, etc.
Meanwhile, for the Berry-Esseen bounds of sample quantiles, one can
refer to Serfling \cite[Theorem 2.3.3 C]{serfling 1980}, Lahiri and
Sun \cite{lahiri 2009}, Yang et al. \cite{yang hu 2012} and the
references therein.

One of the applications of the quantile function is in finance where
many financial returns can be modeled as time series data.
Value-at-risk(VaR) is a popular measure of market risk associated
with an asset or a portfolio of assets. It has been chosen by the
Basel Committee on Banking Supervision as a benchmark risk measure
and has been used by financial institutions for asset management and
minimization of risk. Let $\{X_t\}_{t=1}^n$ be the market value of
an asset over $n$ periods of a time unit, and let $Y_t=\log
(X_t/X_{t-1})$ be the log-returns. Suppose $\{Y_t\}_{t=1}^n$ is a
strictly stationary dependent process with marginal distribution
function $F$. Given a positive value $p$ close to zero, the $1-p$
level VaR is
$$v_p=\inf\{x:F(x)\geq p\},$$
which specifies the smallest amount of loss such that the
probability of the loss in market value being large than $v_p$ is
less than $p$. So, the study of VaR is a well application of $p$th
quantile. Chen and Tang \cite{chen 2005} considered the
nonparametric estimation of VaR and associated standard error
estimation for dependent financial returns. Theoretical properties
of the kernel VaR estimator were investigated in the context of
dependent. For more details, one can refer to Chen and Tang
\cite{chen 2005} and the references therein.

Before stating our works, we need to recall the definition of
$\varphi$-mixing. Let $n$ and $m$ be positive integers. Write
$\mathcal{F}_n^m=\sigma(X_i, n\leq i\leq m)$. Given
$\sigma$-algebras $\mathcal{B},\mathcal{R}$ in $\mathcal{F}$, let $$
\varphi(\mathcal{B},\mathcal{R})=\sup_{A\in \mathcal{B}, B\in
\mathcal{R}, P(A)>0}|P(B|A)-P(B)|.$$ Define the $\varphi$-mixing
coefficients by $$
\varphi(n)=\sup_{k\geq1}\varphi(\mathcal{F}_1^k,\mathcal{F}_{k+n}^\infty),~n\geq0.$$

\textbf{Definition~1.1.} {\it A random variable sequence
~$\{X_n,n\geq1\}$ is said to be a $\varphi$-mixing random variable
sequence if $\varphi(n)\downarrow0$ as $n\rightarrow\infty$.}

In this paper, we go on investigating the Bahadur representation for
sample quantiles under $\varphi$-mixing sequence. Under some mild
conditions such as $\varphi (n)=O(n^{-3})$ and $P(X_i=X_j)=0$,
$\forall~i\neq j$, the rate is established as
$O(n^{-\frac{3}{4}}\log n)$ $a.s.$, which is close to the rate
$O(n^{-\frac{3}{4}}(\log n)^{3/4})$ $a.s.$ in Theorem 1.1 for
$i.i.d.$ random variables. By relaxing the mixing coefficients to
$\sum\nolimits_{n=1}^\infty\varphi^{1/2} (n)<\infty$ and removing
the condition $P(X_i=X_j)=0$, $\forall~i\neq j$, we get the rate as
$O(n^{-1/2}(\log n)^{1/2})$ $a.s.$. The Bahadur representation for
sample quantiles under $\varphi$-mixing sequence has been studied by
Sen \cite{sen 1972}, Babu and Singh \cite{babu 1978} and Yoshihara
\cite{yoshihara 1995}, etc. Comparing our Theorem 2.1 in Section 2
with Theorem 3.1 of Sen \cite{sen 1972}, under some conditions, Sen
obtained the rate $O(n^{-\frac{3}{4}}\log n)$ $a.s.$, where the
mixing coefficients satisfy that $\sum\nolimits_{n=1}^\infty
n^k\varphi^{1/2}(n)<\infty$ for some $k\geq 1$. The mixing
coefficients condition $\varphi(n)=O(n^{-3})$ in our Theorem 2.1 is
relatively weaker.

Under the mixing coefficients condition $\alpha(n)=O(n^{-\beta})$,
$\beta>3$, Zhang et al. \cite{zhang 2013} studied the Bahadur
representation for sample quantiles under $\alpha$-mixing sequence
and obtained the rate $O(n^{-\frac{1}{2}}(\log\log n\cdot\log
n)^{\frac{1}{2}})$, $a.s.$ ( see Theorems 2.3 of Zhang et al. \cite{zhang 2013}). For any $\delta>0$ and $\alpha(n)=O(n^{-\beta})$,
$\beta>\max\{3+\frac{5}{1+\delta},1+\frac{2}{\delta}\}$, they also
obtained the rate
$O(n^{-\frac{3}{4}+\frac{\delta}{4(2+\delta)}}(\log\log n\cdot \log
n)^{\frac{1}{2}})$, $a.s.$ (see Theorem 2.5 of Zhang et al.
\cite{zhang 2013}).

It is a fact that $\varphi$-mixing random variables are
$\alpha$-mixing random variables and $\alpha(n)\leq \varphi(n)$. In
this paper, we investigate the Bahadur representation for sample
quantiles under $\varphi$-mixing sequence. By taking
$\varphi(n)=O(n^{-3})$ in our Theorem 2.2, we get a better rate
$O\left(n^{-\frac{3}{4}}\log n\right)$, $a.s.$ than
$O(n^{-\frac{3}{4}+\frac{\delta}{4(2+\delta)}}(\log\log n\cdot \log
n)^{\frac{1}{2}})$, $a.s.$ obtained by Theorem 2.5 of Zhang et al.
\cite{zhang 2013}. Similarly, by taking $\sum_{n=1}^\infty
\varphi^{1/2}(n)<\infty$ in our Theorem 2.5, we also obtain a better
rate $O(n^{-1/2}(\log n)^{1/2})$, $a.s.$ than
$O(n^{-\frac{1}{2}}(\log\log n\cdot\log n)^{\frac{1}{2}})$, $a.s.$
obtained by Theorems 2.3 of Zhang et al. \cite{zhang 2013}. Although
$\varphi$-mixing random variables are $\alpha$-mixing random
variables, the bounds in our Theorems 2.1-2.5 are better than the
ones obtained by Theorem 2.4, Theorem 2.5 and Theorems 2.1-2.3 of
Zhang et al. \cite{zhang 2013}, respectively.

The organization of this paper is as follows. The main results are
presented in Section 2, some preliminary lemmas are given in Section
3 and the proofs of theorems are provided in Section 4. Throughout
the paper, $C$ and $C_0$ denote positive constants which may be
different in various places. $\lceil x\rceil$ denotes the largest
integer not exceeding $x$. Denote $c_n\sim d_n$, which means
$c_nd_n^{-1}\rightarrow1$ as $n\rightarrow\infty$.

\setcounter{equation}{0}
\section{Main results}
\hspace*{0.5cm} For a fixed $p\in (0,1)$, let $\xi_p=F^{-1}(p)$,
$\xi_{p,n}=F_n^{-1}(p)$.

\textbf{Theorem 2.1.} {\it Let $\{X_n, n\geq1\}$ be a sequence of
$\varphi$-mixing random variables with the mixing coefficients
$\varphi(n)=O(n^{-3})$. Assume that the common marginal distribution
function $F(x)$ possesses a positive continuous density $f(x)$ in a
neighborhood $\mathscr{N}_p$ of $\xi_p$ such that
$0<d=\sup\{f(x):x\in \mathscr{N}_p\}<\infty$. Let $a_n\sim
C_0n^{-\frac{1}{2}}(\log n)^{\frac{3}{4}}$,
$\mathscr{D}_n=[\xi_p-a_n,\xi_p+a_n]$ and
$$H_{p,n}\doteq\sup\limits_{x\in\mathscr{D}_n}|(F_n(x)-F(x))-(F_n(\xi_p)-p)|.$$ Then with
probability 1,
\begin{equation}
H_{p,n}\leq (d+4C_3)n^{-\frac{3}{4}}\log n,~~~\textrm{for
all}~n~\textrm{sufficiently large},\label{n1}
\end{equation}
where $C_3=4[1+4\sum\nolimits_{n=1}^{\infty}\varphi^{1/2}(n)]$. }

\textbf{Theorem 2.2.} {\it Let the conditions of Theorem 2.1 be
satisfied and $P(X_i=X_j)=0$, $\forall~i\neq j$. Assume that
$f^{\prime}(x)$ is bounded in some neighborhood of $\xi_p$, say
$\mathscr{N}_p$. Then with probability 1,
\begin{equation}
\xi_{p,n}=\xi_p-\frac{F_n(\xi_p)-p}{f(\xi_p)}+O\left(n^{-\frac{3}{4}}\log
n\right),~n\rightarrow\infty.\label{n2}
\end{equation} }

\textbf{Theorem 2.3.} {\it  Let $\{X_n, n\geq1\}$ be a sequence of
$\varphi$-mixing random variables with $\sum_{n=1}^\infty
\varphi^{1/2}(n)<\infty$. Assume that the common marginal
distribution function $F(x)$ possesses a positive continuous density
$f(x)$ in a neighborhood $\mathscr{N}_p$ of $\xi_p$ such that
$0<d=\sup\{f(x):x\in \mathscr{N}_p\}<\infty$. For any $\theta>0$,
put $\tau_n=\frac{\sqrt{16C_3+\theta}(\log
n)^{3/2}}{n^{1/2}(\log\log n)^{1/2}}$,
$\mathscr{E}_n=[\xi_p-\tau_n,\xi_p+\tau_n]$, where
$C_3=4[1+4\sum\nolimits_{n=1}^{\infty}\varphi^{1/2}(n)]$. Then with
probability 1
\begin{equation}
\sup\limits_{x\in\mathscr{E}_n}|(F_n(x)-F(x))-(F_n(\xi_p)-p)|\leq
(1+d)\left(\frac{(16C_3+\theta)\log n}{n}\right)^{1/2}\label{n3}
\end{equation} for all $n$ sufficiently
large.}

\textbf{Theorem 2.4.} {\it Let the conditions of Theorem 2.3 hold.
Then with probability 1
\begin{equation}
\sup\limits_{x\in \mathscr{E}_n}|F_n(x)-F(x)|\leq
(1+d)\left(\frac{(16C_3+\theta)\log
n}{4n}\right)^{1/2},~~\textrm{for all $n$ sufficiently
large}.\label{n4}
\end{equation}}

\textbf{Theorem 2.5.} {\it Let the conditions of Theorem 2.3 be
satisfied and $f^{\prime}(x)$ be bounded in some neighborhood of
$\xi_p$. Then with probability 1
\begin{equation}
\xi_{p,n}=\xi_p-\frac{F_n(\xi_p)-p}{f(\xi_p)}+O\left(\frac{(\log
n)^{1/2}}{n^{1/2}}\right),~n\rightarrow\infty.\label{n5}
\end{equation}}

\setcounter{equation}{0}
\section{Preliminary lemmas}

\hspace*{0.5cm} \textbf{Lemma 3.1.}  {\it (Yang et al.
\cite[Corollary A.1]{yang 2012}) Let $\{X_n, n\geq 1\}$ be a
$\varphi$-mixing sequence with $EX_i=0$, $|X_i|\leq d<\infty$, a.s.
$i=1,2,\cdots$, $0<\beta<1$, $m=\lceil n^\beta\rceil$,
$\Delta_2=\sum\nolimits_{i=1}^n EX_i^2$. Then for any
$\varepsilon>0$ and $n\geq 2$, it follows
\begin{equation}
P\left(\left|\sum\limits_{i=1}^nX_i\right|> \varepsilon\right)\leq
2eC_1\exp\left\{-\frac{\varepsilon^2}{2C_2(2\Delta_2+n^\beta
d\varepsilon)}\right\},\nonumber
\end{equation}
where $C_1=\exp\{2en^{1-\beta}\varphi(m)\}$,
$C_2=4[1+4\sum\nolimits_{i=1}^{2m}\varphi^{1/2}(i)]$ and $e$ is the
base of natural logarithm.}

\textbf{Lemma 3.2.} {\it (Serfling \cite[Lemma 1.1.4]{serfling
1980}) Let $F(x)$ be a right-continuous distribution function. The
inverse function $F^{-1}(t)$, $0<t<1$, is nondecreasing and
left-continuous, and satisfies

(i)~~~$F^{-1}(F(x))\leq x,~-\infty<x<\infty$;

(ii)~~$F(F^{-1}(t))\geq t,~0<t<1$;

(iii)~$F(x)\geq t$~\textrm{if and only if}~$x\geq F^{-1}(t)$.}

Inspired by Serfling \cite[Theorem 2.3.2 and Lemma 2.5.4 B]{serfling
1980}, we obtain the following result.

\textbf{Lemma 3.3.} {\it Let $p\in (0,1)$ and $\{X_n, n\geq1\}$ be a
sequence of $\varphi$-mixing random variables with
$\sum_{n=1}^\infty \varphi^{1/2}(n)<\infty$. Assume that the common
marginal distribution function $F(x)$ is differentiable at $\xi_p$,
with $F^{'}(\xi_p)=f(\xi_p)>0$. Suppose that $f^{'}(x)$ is bounded
in a neighborhood of $\xi_p$, say $\mathscr{N}_p$. Then for any
$\delta>0$, with probability 1
\begin{equation}
|\xi_{p,n}-\xi_p|\leq \frac{(2\sqrt C_3+\delta)(\log
n)^{1/2}}{f(\xi_p)n^{1/2}}, ~~\textrm{for all $n$ sufficiently
large},\label{n6}
\end{equation}
where $C_3=4[1+4\sum\nolimits_{n=1}^{\infty}\varphi^{1/2}(n)]$. }

\textbf{Proof.} \ Let $\varepsilon_n=\frac{(2\sqrt C_3+\delta)(\log
n)^{1/2}}{f(\xi_p)n^{1/2}}$, $n>1$, $\frac{1}{3}\leq \beta
<\frac{1}{2}$. Write
$$P(|\xi_{p,n}-\xi_p|>\varepsilon_n)=P(\xi_{p,n}>\xi_p+\varepsilon_n)+P(\xi_{p,n}<\xi_p-\varepsilon_n).$$
By Lemma 3.2 $(iii)$,
\begin{eqnarray*}
P(\xi_{p,n}>\xi_p+\varepsilon_n)&=&P(p>F_n(\xi_p+\varepsilon_n))=P(1-F_n(\xi_p+\varepsilon_n)>1-p)\\
&=&P\left(\sum\limits_{i=1}^nI(X_i>\xi_p+\varepsilon_n)>n(1-p)\right)\\
&=&P\left(\sum\limits_{i=1}^n(V_i-EV_i)>n\delta_{n1}\right),
\end{eqnarray*}
where $V_i=I(X_i>\xi_p+\varepsilon_n)$ and
$\delta_{n1}=F(\xi_p+\varepsilon_n)-p>0$. Likewise, by Lemma 3.2
$(iii)$,
$$
P(\xi_{p,n}<\xi_p-\varepsilon_n)\leq P(p\leq
F_n(\xi_p-\varepsilon_n)) =P\left(\sum\limits_{i=1}^n(W_i-EW_i)\geq
n\delta_{n2}\right),
$$
where $W_i=I(X_i\leq\xi_p-\varepsilon_n)$ and
$\delta_{n2}=p-F(\xi_p-\varepsilon_n)>0$. It is easy to see that
$\{V_i-EV_i, 1\leq i\leq n\}$ and $\{W_i-EW_i, 1\leq i\leq n\}$ are
still $\varphi$-mixing random variables with mean zero and same
mixing coefficients. Since $|V_i-EV_i|\leq 1$,
$\sum\nolimits_{i=1}^n E(V_i-EV_i)^2\leq n$, $|W_i-EW_i|\leq 1$,
$\sum\nolimits_{i=1}^n E(W_i-EW_i)^2\leq n$, it follows from Lemma
3.1 that
$$P(\xi_{p,n}>\xi_p+\varepsilon_n)\leq
2eC_1\exp\left\{-\frac{n\delta_{n1}^2}{2C_3(2+n^\beta\delta_{n1})}\right\},$$
$$
P(\xi_{p,n}<\xi_p-\varepsilon_n)\leq
2eC_1\exp\left\{-\frac{n\delta_{n2}^2}{2C_3(2+n^\beta\delta_{n2})}\right\},
$$
where $C_1=\exp\{2en^{1-\beta}\varphi(m)\}$, $m=\lceil
n^\beta\rceil$ and
$C_3=4[1+4\sum\nolimits_{n=1}^{\infty}\varphi^{1/2}(n)]$.
Consequently,
$$ P(|\xi_{p,n}-\xi_p|>\varepsilon_n)\leq4eC_1
\exp\left\{-\frac{n[\min(\delta_{n1},\delta_{n2})]^2}{2C_3(2+n^\beta[\max(\delta_{n1},\delta_{n2})])}\right\}.
$$ Since $F(x)$ is continuous at $\xi_p$ with $F^{'}(\xi_p)>0$,
$\xi_p$ is the unique solution of $F(x-)\leq p\leq F(x)$ and
$F(\xi_p)=p$. By the assumption on $f^{\prime}(x)$ and Taylor's
expansion,
$$\delta_{n1}=F(\xi_p+\varepsilon_n)-p=f(\xi_p)\varepsilon_n+o(\varepsilon_n)=\frac{(2\sqrt C_3+\delta)(\log
n)^{1/2}}{n^{1/2}}+o(\varepsilon_n),$$ and
$$\delta_{n2}=p-F(\xi_p-\varepsilon_n)=f(\xi_p)\varepsilon_n+o(\varepsilon_n)=\frac{(2\sqrt C_3+\delta)(\log
n)^{1/2}}{n^{1/2}}+o(\varepsilon_n).$$ So
$$
[\min(\delta_{n1},\delta_{n2})]^2\geq(2\sqrt
C_3+\delta/2)^2\cdot\frac{\log n}{n}, ~\textrm{for all} ~n~
\textrm{sufficiently large}.
$$
Since $\frac{1}{3}\leq \beta <\frac{1}{2}$, it has
$n^\beta[\max(\delta_{n1},\delta_{n2})]\rightarrow0$ as
$n\rightarrow\infty$. So we can choose $\delta_1>0$ such that
$$n^\beta[\max(\delta_{n1},\delta_{n2})]\leq \delta_1,~\textrm{and}~\delta_2\doteq\frac{(2\sqrt
C_3+\delta/2)^2}{2C_3(2+\delta_1)}>1$$ for all $n$ sufficiently
large. Hence,
$$P(|\xi_{p,n}-\xi_p|>\varepsilon_n)\leq4eC_1
\exp\left\{-\frac{(2\sqrt C_3+\delta/2)^2}{2C_3(2+\delta_1)}\log
n\right\}=\frac{4eC_1}{n^{\delta_2}}$$
for all $n$ sufficiently
large.

Since $\varphi(n)\downarrow 0$ as $n\rightarrow\infty$ and
$\sum\nolimits_{n=1}^{\infty}\varphi^{1/2}(n)<\infty$, it follows
that $\varphi^{1/2}(n)=o(\frac{1}{n})$. Therefore, $$C_1=\exp\{2e
n^{1-\beta}\varphi(m)\}\leq C \exp\{2e n^{1-\beta}n^{-2\beta}\}\leq
C\exp\{2e\},$$ and
\begin{eqnarray*}
\sum_{n=1}^\infty
P(|\xi_{p,n}-\xi_p|>\varepsilon_n)&=&\sum_{n=1}^{n_0-1}
P(|\xi_{p,n}-\xi_p|>\varepsilon_n)+\sum_{n=n_0}^\infty
P(|\xi_{p,n}-\xi_p|>\varepsilon_n)\\
&\leq& C+4eC_1\sum_{n=n_0}^\infty\frac{1}{n^{\delta_2}}<\infty,
\end{eqnarray*}
which implies that with probability 1 ($wp1$), the
relations $|\xi_{p,n}-\xi_p|>\varepsilon_n$ hold for only finitely
many $n$ by Borel-Cantelli Lemma. Thus (\ref{n6}) holds. $\sharp$

\textbf{Lemma 3.4.} {\it (Wang et al. \cite[Lemma 3.4]{wang 2011})
Let $p\in (0,1)$ and $\xi_{p,n}=F_n^{-1}(p)=\inf\{x: F_n(x)\geq
p\}$. For any $i\neq j$, we assume that $P(X_i=X_j)=0$. Then}
\begin{equation}
p\leq
F_n(\xi_{p,n})<p+\frac{1}{n},~~a.s.\label{a1}
\end{equation}

\setcounter{equation}{0}
\section{Proofs of main results}

\hspace*{0.5cm}\textbf{Proof of Theorem 2.1.} The proof is inspired
by Serfling \cite[Lemma 2.5.4 E]{serfling 1980}. Let $\{b_n,
n\geq1\}$ be a sequence of positive integers such that $b_n\sim
C_0n^{\frac{1}{4}}\log n$ as $n\rightarrow\infty$. For integers
$r=-b_n \cdots, b_n$, put
$$\eta_{r,n}=\xi_p+a_nb_n^{-1}r,~~\alpha_{r,n}=F(\eta_{r+1,n})-F(\eta_{r,n}),$$
and
$$G_{r,n}=F_n(\eta_{r,n})-F_n(\xi_p)-F(\eta_{r,n})+p\doteq\frac{1}{n}\sum_{i=1}^n(Y_i-EY_i),$$
where $Y_i=I(X_i\leq \eta_{r,n})-I(X_i\leq \xi_p)$. Denote
$$g(x)=F_n(x)-F(x)-F_n(\xi_p)+p,~~~x\in \mathscr{D}_n.$$
Then for all $x\in[\eta_{r,n},\eta_{r+1,n}]$,
$$
g(x)\leq F_n(\eta_{r+1,n})- F(\eta_{r,n})-F_n(\xi_p)+p
=G_{r+1,n}+\alpha_{r,n},
$$
and
$$
g(x)\geq F_n(\eta_{r,n})- F(\eta_{r+1,n})-F_n(\xi_p)+p
=G_{r,n}-\alpha_{r,n}.~~$$ So it has
\begin{equation}
H_{p,n}\leq K_n+\beta_n,\label{n7}
\end{equation}
where
$$K_n=\max\{|G_{r,n}|: -b_n\leq r\leq b_n\},~~~~\beta_n=\max\{\alpha_{r,n}: -b_n\leq r\leq b_n-1\}.$$
By the fact
$\eta_{r+1,n}-\eta_{r,n}=a_nb_n^{-1}=n^{-\frac{3}{4}}(\log
n)^{-\frac{1}{4}}$, $-b_n\leq r\leq b_n-1$, we have by the Mean
Value Theorem that
$$\alpha_{r,n}\leq d(\eta_{r+1,n}-\eta_{r,n})=dn^{-\frac{3}{4}}(\log
n)^{-\frac{1}{4}},~-b_n\leq r\leq b_n-1,$$thus
\begin{equation}
0\leq\beta_n\leq dn^{-\frac{3}{4}}(\log n)^{-\frac{1}{4}}.\label{n8}
\end{equation} Taking $\gamma_n=4C_3n^{-\frac{3}{4}}\log n$, we can
check that $\{Y_i-EY_i, 1\leq i\leq n\}$ are still $\varphi$-mixing
random variables with $|Y_i-EY_i|\leq1$, $i=1,2,\cdots,n$. Applying
Lemma 3.1 to $\{Y_i-EY_i, 1\leq i\leq n\}$ and
$\varepsilon=n\gamma_n$, we obtain that
\begin{eqnarray*}
P\left(\left|G_{r,n}\right|>\gamma_n\right)&=&P\left(\left|\sum_{i=1}^n(Y_i-EY_i)\right|>n\gamma_n\right)\\
&\leq&2eC_1\exp\left\{-\frac{n^2\gamma_n^2}{2C_3(2\sum_{i=1}^nEY_i^2+n^\beta\cdot
n\gamma_n)}\right\}\\
&\leq&2eC_1\exp\left\{-\frac{n^2\gamma_n^2}{2C_3(2nz_{r,n}+n^\beta\cdot
n\gamma_n)}\right\}\\
&=&2eC_1\exp\left\{-\frac{n\gamma_n^2}{2C_3(2z_{r,n}+n^\beta\gamma_n)}\right\},
\end{eqnarray*}
where $z_{r,n}=|F(\eta_{r,n})-F(\xi_p)|$. Since
$\varphi(n)=O(n^{-3})$, let $\beta=\frac{1}{4}$ we have that
$C_1=\exp\{2e n^{1-\beta}\varphi(m)\}\leq C\exp\{2e
n^{1-\beta}n^{-3\beta}\}=C\exp\{2e\}$, and
\begin{equation}
P\left(\left|G_{r,n}\right|>\gamma_n\right)\leq2eC_1\exp\left\{-\frac{n\gamma_n^2}{2C_3(2z_{r,n}+n^{\frac{1}{4}}\gamma_n)}\right\}.\label{n9}
\end{equation}
Let $C_4$ be some positive constant such that
$f(\xi_p)<C_4$. Then there exists $N\in \mathcal {N^+}$ such that
$$F(\xi_p+a_n)-F(\xi_p)=f(\xi_p)a_n+o(a_n)<C_4a_n$$
and
$$F(\xi_p)-F(\xi_p-a_n)=f(\xi_p)a_n+o(a_n)<C_4a_n$$
for all $n>N$. Thus
\begin{equation}
z_{r,n}\leq C_4a_n,~\textrm{for}~|r|\leq
b_n~\textrm{and}~n>N.\label{n10}
\end{equation} By (\ref{n9}) and
(\ref{n10}), it follows that
\begin{eqnarray*}
P\left(\left|G_{r,n}\right|>\gamma_n\right)&\leq&2eC_1\exp\left\{-\frac{n\cdot
16C_3^2n^{-\frac{3}{2}}(\log n)^2}{2C_3(2C_4C_0n^{-\frac{1}{2}}(\log
n)^{\frac{3}{4}}+n^{\frac{1}{4}}\cdot 4C_3n^{-\frac{3}{4}}\log
n)}\right\}\\
&\leq&2eC_1\exp\left\{-\frac{
16C_3^2n^{-\frac{1}{2}}(\log n)^2}{9C_3^2n^{-\frac{1}{2}}\log
n}\right\}\\
&=&2eC_1\exp\left\{-\frac{16\log
n}{9}\right\}~=~\frac{2eC_1}{n^{16/9}}
\end{eqnarray*}
for all $n$ sufficiently large. Therefore,
$$
\sum_{n=1}^\infty P\left(K_n>\gamma_n\right)\leq
C+C\sum_{n=n_0}^\infty
\sum_{r=-b_n}^{b_n}P\left(\left|G_{r,n}\right|>\gamma_n\right) \leq
C+C\sum_{n=n_0}^\infty\frac{n^{\frac{1}{4}}\log n}{n^{16/9}}<\infty.
$$
Together with Borel--Cantelli Lemma, it follows that with
probability 1 ($wp1$), the relations $K_n>\gamma_n$ hold for only
finitely many $n$. Hence $wp1$, $K_n\leq\gamma_n$,~~for all $n$
sufficiently large, i.e., $wp1$,
\begin{equation}
 K_n\leq 4C_3n^{-\frac{3}{4}}\log
n,~~~~\textrm{for all}~n~\textrm{sufficiently large}.\label{n11}
\end{equation} Finally, (\ref{n7}), (\ref{n8}) and (\ref{n11})
yield (\ref{n1}). $\sharp$~

\textbf{Proof of Theorem 2.2.} \ By Lemma 3.3, we can see that
$wp1$,
\begin{equation}
 |\xi_{p,n}-\xi_p|\leq \frac{(2\sqrt C_3+\delta)(\log
n)^{1/2}}{f(\xi_p)n^{1/2}},~~\textrm{for all $n$ sufficiently
large},\label{n12}
\end{equation} which implies that $wp1$,
$\xi_{p,n}\in\mathscr{D}_n$, for all $n$ sufficiently large. It
follows from Theorem 2.1 that $wp1$,
\begin{eqnarray*}
|(F_n(\xi_{p,n})-F(\xi_{p,n}))-(F_n(\xi_p)-p)|&\leq&\sup\limits_{x\in\mathscr{D}_n}|(F_n(x)-F(x))-(F_n(\xi_p)-p)|\\
&=&O\left(n^{-\frac{3}{4}}\log n\right),~n\rightarrow\infty,
\end{eqnarray*}
which implies that $wp1$,
\begin{equation}
F_n(\xi_p)-p=F_n(\xi_{p,n})-F(\xi_{p,n})+O\left(n^{-\frac{3}{4}}\log
n\right),~n\rightarrow\infty.\label{n13}
\end{equation} By
(\ref{n12}), (\ref{n13}) and Lemma 3.4, we can obtain that $wp1$,
\begin{eqnarray*}
F_n(\xi_p)-p&=&p+O\left(\frac{1}{n}\right)-F(\xi_{p,n})+O\left(n^{-\frac{3}{4}}\log
n\right)\\&=&-(F(\xi_{p,n})-F(\xi_{p}))+O\left(n^{-\frac{3}{4}}\log
n\right)\\
&=&-f(\xi_{p})(\xi_{p,n}-\xi_{p})-\frac{1}{2}f^{'}(w_n)(\xi_{p,n}-\xi_{p})^2+O\left(n^{-\frac{3}{4}}\log
n\right)\\
&=&-f(\xi_{p})(\xi_{p,n}-\xi_{p})+O\left(n^{-\frac{3}{4}}\log
n\right),~n\rightarrow\infty,
\end{eqnarray*}
where $w_n$ is a random variable between $\xi_{p,n}$ and $\xi_{p}$.
By reorganizing the above equality, $wp1$, (\ref{n2}) holds.
$\sharp$ ~

\textbf{Proof of Theorem 2.3.} For $n>2$, let
$$t_n=\frac{((16C_3+\theta)\log n)^{1/2}}{n^{1/2}},~~\eta_{r,n}=\xi_p+rt_n,~~\Delta_{r,n}=F_n(\eta_{r,n})-F(\eta_{r,n})-F_n(\xi_p)+p$$
for $r=0,\pm1,\pm2,\cdots,\pm\lceil b_n\rceil$ and $b_n=\frac{\log
n}{(\log\log n)^{1/2}}$, $\frac{1}{3}\leq \beta <\frac{1}{2}$.
Denote
$$g(x)=F_n(x)-F(x)-F_n(\xi_p)+p,~~~x\in \mathscr{E}_n.$$
Then for all $x\in[\eta_{r,n},\eta_{r+1,n}]$, it has
\begin{eqnarray*}
g(x)&\leq& F_n(\eta_{r+1,n})- F(\eta_{r,n})-F_n(\xi_p)+p\\
&=&\Delta_{r+1,n}+F(\eta_{r+1,n})-F(\eta_{r,n})~\leq~
\Delta_{r+1,n}+dt_n
\end{eqnarray*}
and
\begin{eqnarray*}
 g(x)&\geq&
F_n(\eta_{r,n})- F(\eta_{r+1,n})-F_n(\xi_p)+p\\
&=&\Delta_{r,n}+F(\eta_{r,n})-F(\eta_{r+1,n})~ \geq~
\Delta_{r,n}-dt_n,~~~~~
\end{eqnarray*}
which imply that
\begin{equation}
\sup\limits_{x\in\mathscr{E}_n}|F_n(x)-F(x)-F_n(\xi_p)+p|\leq\max\limits_{0\leq
|r|\leq \lceil b_n\rceil}|\Delta_{r,n}|+dt_n.\nonumber
\end{equation}
By the notations above, we can see that
\begin{eqnarray*}
P\left(|\Delta_{r,n}|> t_n\right)&=&P\left(|F_n(\eta_{r,n})-F(\eta_{r,n})-F_n(\xi_p)+p|> t_n\right)\\
&\leq& P\left(|F_n(\eta_{r,n})-F(\eta_{r,n})|>
\frac{t_n}{2}\right)+P\left(|F_n(\xi_p)-p|> \frac{t_n}{2}\right)\\
&\doteq&I_n^{(1)}+I_{n}^{(2)}.
\end{eqnarray*}
Denote $\xi_i=I(X_i\leq \eta_{r,n})-EI(X_i\leq
\eta_{r,n}),~i=1,2,\cdots,n.$ Applying Lemma 3.1 to $\xi_1, \xi_2,
\cdots, \xi_n$ and $\varepsilon=nt_n/2$, we have

\begin{eqnarray*}
I_n^{(1)}
&=&P\left(\left|\sum\limits_{i=1}^n\xi_i\right|>\frac{nt_n}{2}
\right)~\leq~2eC_1\exp\left\{-\frac{n^2t_n^2}{8C_3(2n+n^\beta\cdot
nt_n/2)}\right\}\\&\leq&2eC_1\exp\left\{-\frac{nt_n^2}{16C_3+8C_3n^\beta
t_n}\right\}~=~2eC_1\exp\left\{-\frac{(16C_3+\theta)\log
n}{16C_3+8C_3n^\beta t_n}\right\}.
\end{eqnarray*}
Here, $C_1=\exp\{2e n^{1-\beta}\varphi(m)\}\leq C\exp\{2e\}$. It is
easy to check that $n^\beta t_n\rightarrow0$ as
$n\rightarrow\infty$, so we can choose $\theta_1>0$ such that
$$8C_3n^\beta t_n\leq \theta_1<\theta,~\textrm{and}~\theta_2\doteq\frac{16C_3+\theta}{16C_3+\theta_1}>1$$ for all $n$ sufficiently
large. Hence,
\begin{equation}
 I_n^{(1)}\leq \frac{2eC_1}{n^{\theta_2}}~~\textrm{for
all $n$ sufficiently large}.\label{a6}
\end{equation}
Likewise, $I_n^{(2)}=P\left(|F_n(\xi_p)-p|> \frac{t_n}{2}\right)$
satisfies a similar relation. So,
\begin{eqnarray*} \sum\limits_{n=1}^\infty P\left(\max\limits_{0\leq
|r|\leq \lceil b_n\rceil}|\Delta_{r,n}|>t_n\right) &\leq& C+C
\sum\limits_{n=n_0}^\infty \sum\limits_{r=-\lceil
b_n\rceil}^{\lceil b_n\rceil} P\left(|\Delta_{r,n}|> t_n\right)\\
&\leq& C+C\sum\limits_{n=n_0}^\infty \frac{\log n}{(\log \log
n)^{1/2}n^{\theta_2}}<\infty.
\end{eqnarray*}
By Borel-Cantelli Lemma, it follows that $wp1$, the relations
$\max\limits_{0\leq |r|\leq \lceil b_n\rceil}|\Delta_{r,n}|>t_n$
hold for only finitely many $n$. Hence $wp1$
$$
\sup\limits_{x\in\mathscr{E}_n}|F_n(x)-F(x)-F_n(\xi_p)+p|\leq
t_n+dt_n =(1+d)\frac{((16C_3+\theta)\log n)^{1/2}}{n^{1/2}},
$$
for all $n$ sufficiently large.  The proof of (\ref{n3}) is
completed. $\sharp$

\textbf{Proof of Theorem 2.4.}\  For $n>2$, denote
$$t_n=\left(\frac{(16C_3+\theta)\log n)}{4n}\right)^{1/2},~~s_{r,n}=\xi_p+rt_n,~~d_{r,n}=F_n(s_{r,n})-F(s_{r,n}),$$
for $r=0,\pm1,\pm2,\cdots,\pm\lceil b_n\rceil$ and $b_n=\frac{2\log
n}{(\log\log n)^{1/2}}$, $\frac{1}{3}\leq \beta <\frac{1}{2}$. Then
for any $x\in [\xi_p+rt_n,\xi_p+(r+1)t_n]$,
$$d_{r,n}-dt_n\leq F_n(x)-F(x)\leq d_{r+1,n}+dt_n.$$
Hence
\begin{equation}
\sup\limits_{x\in \mathscr{E}_n}|F_n(x)-F(x)|\leq
\max\limits_{0\leq|r|\leq \lceil
b_n\rceil}|d_{r,n}|+dt_n.\label{n14}
\end{equation}
Denote $\eta_i=I(X_i\leq \xi_p+rt_n)-EI(X_i\leq \xi_p+rt_n)$, $1\leq
i\leq n$. Similar to the proof of (\ref{a6}), by Lemma 3.1, we have
\begin{eqnarray*}
P(|d_{r,n}|>t_n)&=&P\left(\left|F_n(\xi_p+rt_n)-F(\xi_p+rt_n)\right|>t_n\right)~=~P\left(\left|\sum_{i=1}^n \eta_i\right|>nt_n\right)\\
&\leq& 2eC_1\exp\left\{-\frac{(16C_3+\theta)\log
n}{16C_3+8C_3n^\beta t_n}\right\} \leq
\frac{2eC_1}{n^{\theta_2}},~~~\textrm{for all $n$ sufficiently
large},
\end{eqnarray*}
where $C_1=\exp\{2e n^{1-\beta}\varphi(m)\}<\infty$ and
$\theta_2\doteq\frac{16C_3+\theta}{16C_3+\theta_1}>1$.
 Therefore,
\begin{equation*}
\sum\limits_{n=1}^\infty P\left(\max\limits_{0\leq|r|\leq \lceil
b_n\rceil}|d_{r,n}|>t_n\right)\leq C+C \sum\limits_{n=n_0}^\infty
\frac{\log n}{(\log \log n)^{1/2}n^{\theta_2}}<\infty.
\end{equation*}
By Borel-Cantelli Lemma, we obtain that $wp1$, the relations
$\max\limits_{0\leq |r|\leq \lceil b_n\rceil}|d_{r,n}|>t_n$ hold for
only finitely many $n$. Together with (\ref{n14}), we can get
(\ref{n4}) immediately. $\sharp$

\textbf{Proof of Theorem 2.5.} Lemma 3.3 implies that $wp1$,
\begin{equation}
|\xi_{p,n}-\xi_p|\leq \frac{(2\sqrt C_3+\delta)(\log
n)^{1/2}}{f(\xi_p)n^{1/2}}<\tau_n, ~~\textrm{for all $n$
sufficiently large},\label{n15}
\end{equation}
$i.e.$ $wp1$, $\xi_{p,n}\in \mathscr{E}_n$ for all $n$ sufficiently
large. It follows from Theorem 2.3 that $wp1$,
\begin{equation}
F_n(\xi_p)-p=F_n(\xi_{p,n})-F(\xi_{p,n})+O\left(\frac{(\log
n)^{1/2}}{n^{1/2}}\right),~~n\rightarrow\infty.\label{n16}
\end{equation}
By (\ref{n15}), assumption on $f^{\prime}(x)$, Taylor's expansion
and Theorem 2.4, we have that $wp1$,
\begin{eqnarray*}
|F_n(\xi_{p,n})-F(\xi_p)|&\leq&|F_n(\xi_{p,n})-F(\xi_{p,n})|+|F(\xi_{p,n})-F(\xi_p)|\\
&\leq& \sup\limits_{x\in
\mathscr{E}_n}|F_n(x)-F(x)|+f(\xi_p)|\xi_{p,n}-\xi_p|+o(|\xi_{p,n}-\xi_p|)\\
&=&O\left(\frac{(\log
n)^{1/2}}{n^{1/2}}\right),~~n\rightarrow\infty.
\end{eqnarray*}
Together with (\ref{n16}), we obtain that $wp1$,
\begin{eqnarray*}
F_n(\xi_p)-p&=&F(\xi_p)-F(\xi_{p,n})+O\left(\frac{(\log n)^{1/2}}{n^{1/2}}\right)\\
&=&-f(\xi_p)(\xi_{p,n}-\xi_p)-\frac{1}{2}f^{\prime}(w_n)(\xi_{p,n}-\xi_p)^2+O\left(\frac{(\log
n)^{1/2}}{n^{1/2}}\right)\\
&=&-f(\xi_p)(\xi_{p,n}-\xi_p)+O\left(\frac{(\log
n)^{1/2}}{n^{1/2}}\right),~~n\rightarrow\infty,
\end{eqnarray*}
where
$w_n$ is a random variable between $\xi_{p,n}$ and $\xi_p$.
Reorganizing the above equality, we obtain that $wp1$,
$$\xi_{p,n}-\xi_p=-\frac{F_n(\xi_p)-p}{f(\xi_p)}+O\left(\frac{(\log
n)^{1/2}}{n^{1/2}}\right)~~n\rightarrow\infty. $$ The proof of
(\ref{n5}) is completed. $\sharp$

\end{document}